\documentclass[a4paper, oneside, 11pt]{article}%
\usepackage[pagewise]{lineno}
\usepackage{amsmath}
\usepackage{amsfonts}
\usepackage{amssymb}
\usepackage{graphicx}
\usepackage{appendix}
\usepackage[square, numbers, sort&compress]{natbib}
\usepackage{color}%
\providecommand{\U}[1]{\protect\rule{.1in}{.1in}}

\newtheorem{theorem}{Theorem}[section]

\newtheorem{corollary}[theorem]{Corollary}

\newtheorem{assumption}{Assumption}

\newtheorem{lemma}[theorem]{Lemma}

\newtheorem{remark}[theorem]{Remark}

\numberwithin{equation}{section}

\newcommand{\E}{{\mathbb E}}
\newcommand{\R}{{\mathbb R}}

\newcommand{\eof}{\hfill{$\Box$}}

\newcommand{\BMO}{L^{2, \;\mathrm{BMO}}_{\mathcal F^{W}}(0, T;\mathbb{R}^{n})}
\usepackage{url}
\usepackage{hyperref}
\hypersetup{
colorlinks=true, 
breaklinks=true, 
urlcolor=green , 
linkcolor=red, 
citecolor=blue, 
pdftitle={}, 
pdfauthor={}, 
pdfsubject={} 
}
\newcommand{\esssup}{\ensuremath{\operatorname*{ess\:sup}}}

\newcommand{\cM}{\ensuremath{\mathcal{M}}}

\newcommand{\ass}{\gamma}

\usepackage[normalem]{ulem}

\begin{document}

\title{Non-homogeneous stochastic LQ control with regime switching and random coefficients}
\author{Ying Hu \thanks{Univ Rennes, CNRS, IRMAR-UMR 6625, F-35000 Rennes, France. Email:
\texttt{ying.hu@univ-rennes1.fr }}
\and Xiaomin Shi\thanks{Corresponding author. School of Statistics and Mathematics, Shandong University of Finance and Economics, Jinan
250100, China. Email: \texttt{shixm@mail.sdu.edu.cn}}
\and Zuo Quan Xu\thanks{Department of Applied Mathematics, The Hong Kong Polytechnic University, Kowloon, Hong Kong.
Email: \texttt{maxu@polyu.edu.hk}}}

\maketitle
This paper is concerned with a general non-homogeneous stochastic linear quadratic (LQ) control problem with regime switching and random coefficients. We obtain the explicit optimal state feedback control and optimal value for this problem in terms of two systems of backward stochastic differential equations (BSDEs): one is the famous stochastic Riccati equation and the other one is a new linear multi-dimensional BSDE with all coefficients being unbounded. The existence and uniqueness of the solutions to these two systems of BSDEs are proved by means of BMO martingales and contraction mapping method. At last, the theory is applied to study an asset-liability management problem under the mean-variance criteria.

{\textbf{Keywords}.} Non-homogeneous stochastic LQ problem, regime switching, BSDE, unbounded coefficients, mean-variance, asset-liability management.

\textbf{Mathematics Subject Classification (2020)} 93E20 60H30 91G10

\addcontentsline{toc}{section}{\hspace*{1.8em}Abstract}

\section{Introduction}

Since the pioneering work of Wonham \cite{Wo}, stochastic linear-quadratic (LQ) theory has been extensively studied by numerous researchers. For instance, Bismut \cite{Bi} was the first one who studied stochastic LQ problems with random coefficients. In order to obtain the optimal random feedback control, he formally derived a stochastic Riccati equation (SRE). But he could not solve the SRE in the general case. It is Kohlmann and Tang \cite{KT}, for the first time, that established the existence and uniqueness of the one-dimensional SRE.
Tang \cite{Ta} made another breakthrough and proved the existence and uniqueness of the matrix valued SRE with uniformly positive control weighting matrix. Chen, Li and Zhou \cite{CLZ}, Sun, Xiong and Yong \cite{SXY} studied the indefinite stochastic LQ problem which is different obviously from its deterministic counterpart. Kohlmann and Zhou \cite{KZ} established the relationship between stochastic LQ problems and backward stochastic differential equations (BSDEs). Hu and Zhou \cite{HZ} solved the stochastic LQ problem with cone control constraint. Please refer to Chapter 6 in Yong and Zhou \cite{YZ} for a systematic accounts on this subject.

Stochastic LQ problems for Markovian regime switching system were studied in
Li and Zhou \cite{LxZ}, Wen, Li and Xiong \cite{WLX} and Zhang, Li and Xiong \cite{ZLX} where sufficient and necessary conditions of the existence of optimal control, weak closed-loop solvability, open-loop solvability and closed-loop solvability were established. But the coefficients are assumed to be \emph{deterministic} functions of time $t$ for each given regime $i$ in the above three papers,
so their Riccati equation systems are indeed deterministic ordinary differential equations (ODEs).
Our previous work \cite{HSX} studied a cone-constrained stochastic LQ problem with regime switching in which the coefficients are \emph{stochastic} processes for each give regime $i$.
Due to the randomness of the coefficients, we have to solve stochastic Riccati equations which are actually a new type of BSDEs.

This paper further explores general stochastic LQ problem with regime switching and random coefficients. Compared with our previous work \cite{HSX}, non-homogeneous terms emerge in both the state process and cost functional in the present LQ problem. Two related systems of (\emph{multi-dimensional}) BSDEs are introduced: the first one is the so called system of SRE whose solvability is established by slightly modifying our previous argument in \cite{HSX}. By contrast, the existing argument cannot deal with the second one because its coefficients, which depend on the solution of the first one, are \emph{inevitably unbounded}.
A similar linear BSDE with unbounded coefficients has appeared in \cite{HSX}, its solvability is established by change of variables so that the new one becomes a linear BSDE with partially bounded coefficients that can be dealt by contraction mapping.
For our problem, all the coefficients of the linear BSDE are unbounded (see Remarks \ref{oncoefficients1} and \ref{oncoefficients2}) and we cannot make change of variables to reduce the BSDE to a solvable one, hence new method is called for solving it.
The main idea to establish the solvability of the new type of BSDEs is first to get some estimates of BMO martingales and then establish the result for one-dimensional system,
and finally apply contraction mapping method to get the result for multi-dimensional system. This constitutes the major technique contribution of this paper.
Eventually we obtain the optimal feedback control and optimal value of the LQ problem through these two systems of BSDEs and some verification arguments.
We establish the above results in both a standard case and a singular case.

On the other hand, asset-liability management (ALM) is one important class of problems in risk management and insurance.
We will apply our theory to such a problem under mean-variance criteria.
Continuous time mean-variance problems have been extensively studied by LQ optimal control theory; see, e.g. in \cite{BJPZ}, \cite{JYZ}, \cite{LZL}, \cite{LZ}, \cite{ZLi}, \cite{ZY} and the references therein. Chiu and Li \cite{CL} investigated firstly in a continuous time setting an ALM problem under the mean-variance criteria. Xie, Li and Wang \cite{XLW} considered this problem with liability process driven by another correlated Brownian motion.
Wei and Wang \cite{WW} found a time-consistent open-loop equilibrium strategy for the problem.
Zeng and Li \cite{ZL} studied this problem in a jump diffusion market.
Chen, Yang and Yin \cite{CYY} generalized the model of \cite{CL} to a setting where the coefficients and liability process were modulated by a continuous time Markov chain and geometric Brownian motion respectively. With liability being described by drifted Brownian motion, Xie \cite{Xie} studied a mean-variance ALM problem with deterministic and Markov chain modeled coefficients.

In the above Markov chain modulated models, the market parameters, such as the interest rate, stock appreciation rates and volatilities are assumed to be \emph{deterministic} functions of time $t$ for each given regime $i$.
Again, their Riccati equation systems are indeed deterministic ODEs.
In practice, however, these market parameters are affected by the uncertainties caused by noises. Thus, it is too restrictive to set market parameters as deterministic even if the market status is known. From practical point of view, it is necessary to allow the market parameters to depend on both the noises and the Markov chain. Shen, Wei and Zhao \cite{SWZ} studied a mean-variance ALM problem under non-Markovian regime switching model. They characterized the optimal portfolios in terms of two \emph{one-dimensional} BSDEs with jumps. In this paper,
we study a mean-variance ALM problem with regime switching and random coefficients using the non-homogeneous stochastic LQ control theory that will be established. The optimal portfolios are characterized by two systems of (\emph{multi-dimensional}) BSDEs without jumps.

The rest part of this paper is organized as follows. In Section \ref{Problem formulation}, we introduce the non-homogeneous stochastic LQ problem with regime switching and random coefficients. Section \ref{unbounded coefficients} is devoted to establishing the global solvability of two systems of BSDEs arisen from the above LQ problem. Section \ref{sec:solution} provides the solution to the original stochastic LQ problem. In Section \ref{sec:application}, we apply the general results to solve a mean-variance ALM problem.
\section{Problem formulation}\label{Problem formulation}
Let $(\Omega, \mathcal F, \mathbb{P})$ be a fixed complete probability space on which are defined a standard $n$-dimensional Brownian motion $W(t)=(W_1(t), \ldots,$ $W_n(t))'$ and a continuous-time stationary Markov chain $\alpha_t$ valued in a finite state space $\mathcal M=\{1, 2, \ldots, \ell\}$ with $\ell>1$. We assume that $\{W(t)\}_{t\geq0}$ and $\{\alpha_t\}_{t\geq0}$ are independent processes. The Markov chain has a generator $Q=(q_{ij})_{\ell\times \ell}$ with $q_{ij}\geq 0$ for $i\neq j$ and $\sum_{j=1}^{\ell}q_{ij}=0$ for every $i\in\mathcal{M}$.
Define the filtrations $\mathcal F_t=\sigma\{W(s), \alpha_s: 0\leq s\leq t\}\bigvee\mathcal{N}$ and $\mathcal F^W_t=\sigma\{W(s): 0\leq s\leq t\}\bigvee\mathcal{N}$, where $\mathcal{N}$ is the totality of all the $\mathbb{P}$-null sets of $\mathcal{F}$.

Throughout this paper, we denote by $\R^n$ the set of $n$-dimensional column vectors, by $\R^n_+$ the set of vectors in $\R^n$ whose components are nonnegative, by $\R^{m\times n}$ the set of $m\times n$ real matrices, and by $\mathbb{S}^n$ the set of symmetric $n\times n$ real matrices. For $x\in\R$, we define $x^+:=\max\{x, 0\}$, and $x^-:=\max\{-x, 0\}$. If $M=(m_{ij})\in \R^{m\times n}$, we denote its transpose by $M'$, and its norm by $|M|=\sqrt{\sum_{ij}m_{ij}^2}$. If $M\in\mathbb{S}^n$ is positive definite (positive semidefinite) , we write $M>$ ($\geq$) $0.$ We write $A>$ ($\geq$) $B$ if $A, B\in\mathbb{S}^n$ and $A-B>$ ($\geq$) $0.$
We will use the following notations throughout the paper:
\begin{align*}
L^{\infty}_{\mathcal{F}}(\Omega;\mathbb{R})&=\Big\{\xi:\Omega\rightarrow
\mathbb{R}\;\Big|\;\xi\mbox { is }\mathcal{F}_{T}\mbox{-measurable, and essentially bounded}\Big\}, \\
L^{2}_{\mathcal F}(0, T;\mathbb{R})&=\Big\{\phi:[0, T]\times\Omega\rightarrow
\mathbb{R}\;\Big|\;\phi(\cdot)\mbox{ is an }\{\mathcal{F}%
_{t}\}_{t\geq0}\mbox{-adapted process }\\
&\qquad\mbox{with }\E\int_{0}^{T}|\phi(t)|^{2}dt<\infty
\Big\}, \\
L^{\infty}_{\mathcal{F}}(0, T;\mathbb{R})&=\Big\{\phi:[0, T]\times\Omega
\rightarrow\mathbb{R}\;\Big|\;\phi(\cdot)\mbox{ is an }\{\mathcal{F}%
_{t}\}_{t\geq0}\mbox{-adapted essentially}\\
&\qquad\mbox{bounded process} \Big\}.
\end{align*}

These definitions are generalized in the obvious way to the cases that $\mathcal{F}$ is replaced by $\mathcal{F}^W$ and $\mathbb{R}$ by $\mathbb{R}^n$, $\mathbb{R}^{n\times m}$ or $\mathbb{S}^n$.
In our argument, $t$, $\omega$, ``almost surely" (a.s.) and ``almost everywhere" (a.e.) may be suppressed for notation simplicity in some circumstances when no confusion occurs.


Consider the following $\mathbb{R}$-valued linear stochastic differential equation (SDE):
\begin{align}
\label{state}
\begin{cases}
dX(t)=\left[A(t, \alpha_t)X(t)+B(t, \alpha_t)'u(t)+b(t,\alpha_t)\right]dt\\
\qquad\qquad+\left[C(t, \alpha_t)'X(t)+u(t)'D(t, \alpha_t)'+\rho(t,\alpha_t)'\right]dW(t), \ t\geq0, \\
X(0)=x, \ \alpha_0=i_0,
\end{cases}
\end{align}
where $A(t, \omega, i), \ B(t, \omega, i), \ b(t,\omega,i), \ C(t, \omega, i), \ D(t, \omega, i), \ \rho(t,\omega,i)$ are all $\{\mathcal{F}^W_t\}_{t\geq 0}$-adapted processes of suitable sizes for $i\in\cM$, the initial states $x\in\mathbb{R}$ and $i_{0}\in\cM$ are known.

The class of admissible controls is defined as the set
\begin{align*}
\mathcal{U}:= L^2_\mathcal{F}(0, T;\mathbb{R}^m).
\end{align*}

If $u(\cdot)\in\mathcal{U}$ and $X(\cdot)$ is the associated (unique strong) solution of \eqref{state}, then we refer to $(X(\cdot), u(\cdot))$ as an admissible pair.

The general stochastic linear quadratic optimal control problem (stochastic LQ problem, for short) is stated as follows:
\begin{align}
\begin{cases}
\mathrm{Minimize} &\ J(x, i_0, u(\cdot))\\
\mbox{subject to} &\ (X(\cdot), u(\cdot)) \mbox{\ is admissible for} \ \eqref{state},
\end{cases}
\label{LQ}%
\end{align}
where the cost functional is given as the following quadratic form
\begin{align}\label{costfunctional}
J(x, i_0, u(\cdot)):=&\mathbb{E}\Big[\int_0^T\Big(Q(t, \alpha_t)(X(t)-q(t,\alpha_t))^2\nonumber\\
&\quad+(u(t)-p(t,\alpha_t))'R(t, \alpha_t)(u(t)-p(t,\alpha_t))\Big)dt\nonumber\\
&\quad+G(\alpha_T)\big(X(T)-g(\alpha_T)\big)^2\Big].
\end{align}

The associated value function is defined as
\begin{align*}
V(x,i_0):=\inf_{u\in\mathcal{U}}J(x, i_0, u(\cdot)), \ x\in\mathbb{R}, \ i_0\in\cM.
\end{align*}

To make sure the well-posedness of the LQ problem \eqref{LQ}, we put the following assumptions.
\begin{assumption} \label{assu1}
\emph{For all $i\in\cM$,}
{\small\begin{align*}
\begin{cases}
A(t, \omega, i), \ b(t, \omega, i), \ q(t, \omega, i)\in L_{\mathcal{F}^W}^\infty(0, T;\mathbb{R}), \ B(t, \omega, i), \ p(t, \omega, i)\in L_{\mathcal{F}^W}^\infty(0, T;\mathbb{R}^m), \\
C(t, \omega, i), \ \rho(t, \omega, i)\in L_{\mathcal{F}^W}^\infty(0, T;\mathbb{R}^n), \ D(t, \omega, i)\in L_{\mathcal{F}^W}^\infty(0, T;\mathbb{R}^{n\times m}), \\
Q(t, \omega, i)\in L_{\mathcal{F}^W}^\infty(0, T;\mathbb{R}), \ R(t, \omega, i)\in L_{\mathcal{F}^W}^\infty(0, T;\mathbb{S}^m),\\
G(\omega, i)\in L_{\mathcal{F}^W}^\infty(\Omega;\mathbb{R}),\ g(\omega,i)\in L_{\mathcal{F}^W}^\infty(\Omega;\mathbb{R}).
\end{cases}
\end{align*}}
\end{assumption}

\begin{assumption}\label{assu2}
\emph{There exists a constant $\delta>0$ such that at least one of the following cases holds.}
\begin{itemize}
\item[(i)] \emph{\textbf{Standard case.} $Q(i)\geq0$, $R(i)\geq \delta I_{m}$ and $G(i)\geq\delta$.}
\item[(ii)] \emph{\textbf{Singular case.} $Q(i)\geq0$, $ R(i)\geq0$, $G(i)\geq\delta$ and $D(i)'D(i)\geq\delta I_{m}$ for all $i\in\cM$, where $I_{m}$ denotes the $m$-dimensional identity matrix.} 
\end{itemize}
\end{assumption}

Under Assumption \ref{assu2}, clearly we have $J(x, i_0, u(\cdot))\geq 0$, for all $(x, i_0, u)\in\R\times\cM\times\mathcal{U}$.
The LQ problem \eqref{LQ} is said to be solvable, if there exists a control $u^*(\cdot)\in\mathcal{U}$ such that
\begin{align*}
J(x, i_0, u^*(\cdot))\leq J(x, i_0, u(\cdot)), \quad \forall\; u(\cdot)\in\mathcal{U},
\end{align*}
in which case, $u^*(\cdot)$ is called an optimal control for the LQ problem \eqref{LQ}, and the optimal value is
\begin{align*}
V(x, i_0)=J(x, i_0, u^*(\cdot)).
\end{align*}
\section{Linear BSDEs with unbounded coefficients}\label{unbounded coefficients}
To tackle the LQ problem \eqref{LQ}, we first introduce the following system of ($\ell$-dimensional) BSDEs (remind that the arguments $t$ and $\omega$ are suppressed):
\begin{align}
\label{P1}
\begin{cases}
dP(i)=-\Big[(2A(i)+C(i)'C(i))P(i)+2C(i)'\Lambda(i)+Q(i)\\
\quad\quad\quad\quad+H(P(i), \Lambda(i), i)+\sum\limits_{j=1}^\ell q_{ij}P(j)\Big]dt+\Lambda(i)'dW, \\
P(T,i)=G(i),\\
R(i)+P(i)D(i)'D(i)>0, \ \mbox{ for all $i\in\cM$},
\end{cases}
\end{align}

\noindent where
{\small\begin{align*}
&\quad H(t, \omega, P, \Lambda, i)\\
&=-\big(PB(i)+D(i)'(PC(i)+\Lambda)\big)'\big(R(i)+PD(i)'D(i)\big)^{-1}\big(PB(i)+D(i)'(PC(i)+\Lambda)\big).
\end{align*}}

The equation \eqref{P1} is referred to as the stochastic Riccati equation for the LQ problem \eqref{LQ}. By a solution to \eqref{P1}, we mean a $2\ell$-dimensional adapted processes $(P(i),\Lambda(i))_{i=1}^\ell$ satisfying \eqref{P1} and $(P(i),\Lambda(i))\in L^\infty_{\mathcal{F}^W}(0,T; \mathbb{R})\times L^{2}_{\mathcal {F}^W}(0, T;\mathbb{R}^n)$ for all $i\in\cM$. Furthermore, a solution of \eqref{P1} is called nonnegative (resp. uniformly positive) if $P(i)\geq 0$ (resp. $P(i)\geq c$\footnote{We shall use $c$ to represent a generic positive constant which does not depend on $t$ or $i$ and can be different from line to line.} for some constant $c>0$) for all $i\in\cM$.

To show the above BSDE has a solution in the sequel, we need the concept of BMO martingales.
Here we recall some facts about BMO martingales; see Kazamaki \cite{Ka}. A process $\int_0^\cdot \Lambda(s)'dW(s)$ is a BMO martingale on $[0,T]$ if and only if its $\mathrm{BMO}_2$ normal on $[0,T]$ is finite, namely,
\[\left\Vert\int_0^\cdot \Lambda(s)'dW(s)\right\Vert_{\mathrm{BMO}_2}:=\sup_{\tau\leq T}\left(\esssup\mathbb{E}\bigg[\int_\tau^T|\Lambda(s)|^2ds\;\Big|\;\mathcal F^W_\tau\bigg]\right)^{\frac{1}{2}}<\infty,\]
here and hereafter the $\sup_{\tau\leq T}$ is taken over all $\{\mathcal{F}^W_t\}_{t\geq 0}$-stopping times $\tau\leq T$.
The Dol\'eans-Dade stochastic exponential $\mathcal E(\int_0^\cdot \Lambda(s)'dW(s))$ of a BMO martingale $\int_0^\cdot \Lambda(s)'dW(s)$ is a uniformly integrable martingale. Moreover, if $\int_0^\cdot \Lambda(s)'dW(s)$ and $\int_0^\cdot Z(s)'dW(s)$ are both BMO martingales, then under the probability measure $\widetilde{\mathbb{P}}$ defined by $\frac{d\widetilde{\mathbb{P}}}{d\mathbb{P}}\big|_{\mathcal{F}^W_T}=\mathcal E \big(\int_0^T Z(s)'dW(s)\big)$, $\widetilde W(\cdot):=W(\cdot)-\int_0^\cdot Z(s)ds$ is a standard Brownian motion, and $\int_0^\cdot \Lambda(s)'d\widetilde W(s)$ is a BMO martingale.
The following space plays an important role in our argument
{\small\begin{align*}
\BMO&:=\bigg\{\phi \in L^{2}_{\mathcal F^{W}}(0, T;\mathbb{R}^{n}) \;\bigg|\; \int_0^\cdot\phi(s)'dW(s) \mbox{ is a BMO martingale}\bigg\}.
\end{align*}}
\begin{lemma}
\label{Riccatisingular}
Under Assumptions \ref{assu1} and \ref{assu2}, the system of BSDEs \eqref{P1} admits a unique uniformly positive solution $(P(i), \Lambda(i))_{i\in\cM}$.
\end{lemma}

\noindent \emph{Proof.} According to Theorems 3.5 (resp. Theorem 3.6) of \cite{HSX}, there exists a unique nonnegative (resp. uniformly positive) solution $(P(i), \Lambda(i))_{i=1}^{\ell}$ to BSDE \eqref{P1} under Assumptions \ref{assu1} and \ref{assu2} (i) (resp. \ref{assu2} (ii)). Note that Assumption \ref{assu2} (i) is stronger than the standard assumption in Theorem 3.5 of \cite{HSX}.
So it remains to show that the solution of \eqref{P1} is actually uniformly positive under Assumptions \ref{assu1} and \ref{assu2} (i).

Let $c_1$, $c_2$ be two positive constants such that $P(i)\leq c_1$, and
\[
2A(i)+C(i)'C(i)+q_{ii}-\frac{2c_1}{\delta}|B(i)+D(i)'C(i)|^2>-c_2, \ \mbox{for all} \ i\in\cM.
\]

Consider the following $\ell$-dimensional BSDE:
\begin{align}
\label{P1Ndecouple}
\begin{cases}
d\underline P(i)=-\Big[(2A(i)+C(i)'C(i)+q_{ii})\underline P(i)+2C(i)'\underline \Lambda(i)+Q(i)\\
\qquad\qquad+H(\underline P(i),\underline \Lambda(i), i)\Big]dt+\underline \Lambda(i)'dW, \\
\underline P(T,i)=\delta,\\
R(i)+\underline P(i)D(i)'D(i)>0, \ \mbox{ for all} \ i\in\cM.
\end{cases}
\end{align}

This is a decoupled system of BSDEs.
From Theorem 4.1 and Theorem 5.2 of \cite{HZ}, the $i$th equation in \eqref{P1Ndecouple} admits a unique, hence maximal solution (see page 565 of \cite{Ko} for its definition) $(\underline P(i), \ \underline\Lambda(i))\in L^\infty_{\mathcal{F}^W}(0, T; \mathbb {R})\times \BMO$, and $\underline P(i)\geq0$ for all $i\in\cM$. From the proof of Theorem 3.5 of \cite{HSX}, the solution $(P(i), \ \Lambda(i))_{i\in\cM}$ of \eqref{P1} could be approximated by solutions of a sequence of BSDEs with Lipschitz generators\footnote{As for a BSDE $Y(t)=\xi+\int_t^T f(s, Y, Z)ds-\int_t^T Z'dW(s)$, the function $f$ is called the generator and the random variable $\xi$ is called the terminal value.}. Thus we can use comparison theorem for multi-dimensional BSDEs (see e.g. Lemma 3.4 of \cite{HSX}) and then pass to the limit to get
\begin{align}\label{lowerbound}
P(i)\geq\underline P(i), \ \mbox{for all} \ i\in\cM.
\end{align}

Let $g:\mathbb{R}^+\rightarrow [0, 1]$ be a smooth truncation function satisfying $g(x)=1$ for $x\in[0, c_1]$, and $g(x)=0$ for $x\in[2c_1, +\infty)$. Notice that
$c_1\geq P(i)\geq\underline P(i)\geq 0$, so $(\underline P(i), \ \underline\Lambda(i))$ is still a solution of the $i$th equation in BSDE \eqref{P1Ndecouple} with $H(P,\Lambda,i)$ replaced by $H(P,\Lambda,i)g(P)$ in the generator.

Notice that for $P=\underline P(i), \ \Lambda=\underline\Lambda(i)$, we have, under Assumptions \ref{assu1} and \ref{assu2} (i),
\begin{align*}
&\quad\;H(P, \Lambda, i)g(P)\\
&\geq-\frac{1}{\delta}|PB(i)+PD(i)'C(i)+D(i)'\Lambda|^2g(P)\\
&=-\frac{P^2}{\delta}|B(i)+D(i)'C(i)|^2g(P)
-\frac{2P}{\delta}\big(B(i)+D(i)'C(i)\big)'D(i)'\Lambda g(P)\\
&\quad-\frac{1}{\delta}|D(i)'\Lambda|^2 g(P)\\
&\geq-\frac{2c_1P}{\delta}|B(i)+D(i)'C(i)|^2
\\
&\quad-\frac{2P}{\delta}\big(B(i)+D(i)'C(i)\big)'D(i)'\Lambda g(P)-\frac{1}{\delta}|D(i)'\Lambda|^2 g(P).
\end{align*}

The following BSDE
\begin{align*}
\begin{cases}
\quad dP=-\Big[-c_2P+2C(i)'\Lambda
-\frac{2P}{\delta}\big(B(i)+D(i)'C(i)\big)'D(i)'\Lambda g(P)\\
\qquad\quad\quad-\frac{1}{\delta}|D(i)'\Lambda|^2 g(P)\Big]dt+\Lambda'dW,\\
P(T)=\delta,
\end{cases}
\end{align*}
has a Lipschitz generator,
thus it admits a unique solution $(\delta e^{-c_2(T-t)},0)$. Then the maximal solution argument (Theorem 2.3 of \cite{Ko}) gives
\[
\underline P(t,i)\geq\delta e^{-c_2(T-t)}\geq\delta e^{-c_2T}.
\]

Combining with \eqref{lowerbound}, we proved that the solution $(P(i), \Lambda(i))_{i\in\cM}$ of \eqref{P1} is actually uniformly positive under Assumptions \ref{assu1} and \ref{assu2} (i).
\eof

\bigskip

In addition to the stochastic Riccati equation \eqref{P1},
we need to consider another system of BSDEs in order to solve the non-homogeneous stochastic LQ problem \eqref{LQ}.

Let $(P(i),\Lambda(i))_{i\in\cM}$ be the unique uniformly positive solution to \eqref{P1}.
Set
\[
\Gamma(i)=\big(R(i)+P(i)D(i)'D(i)\big)^{-1}\big(P(i)B(i)+D(i)'(P(i)C(i)+\Lambda(i))\big).
\]

We consider the following system of ($\ell$-dimensional) linear BSDEs,
\begin{align}
\label{K}
\begin{cases}
dK(i)=-\Big[\big(P(i)B(i)+D(i)'(P(i)C(i)+\Lambda(i))\big)'\big(R(i)+P(i)D(i)'D(i)\big)^{-1}\\
\qquad\qquad\times\big[D(i)'(P(i)\rho(i)-L(i))-K(i)B(i)-R(i)p(i)\big]\\
\qquad\qquad +A(i)K(i)+C(i)'L(i)-P(i)(C(i)'\rho(i)+b(i))-\rho(i)'\Lambda(i)\\
\qquad\qquad+q(i)Q(i)+\sum\limits_{j=1}^\ell q_{ij}K(j)\Big]dt +L(i)'dW\\
\qquad\quad=-\Big[\big(A(i)-B(i)'\Gamma(i)\big)K(i)+\big(C(i)-D(i)\Gamma(i)\big)'L(i)+(P(i)D(i)'\rho(i)\\
\qquad\qquad-R(i)p(i))'\Gamma(i)+q(i)Q(i)-P(i)(C(i)'\rho(i)+b(i))\\
\qquad\qquad-\rho(i)'\Lambda(i)+\sum\limits_{j=1}^\ell q_{ij}K(j)\Big]dt+L(i)'dW,\\
K(T,i)=G(i)g(i),\ \mbox{for all} \ i\in\cM.
\end{cases}
\end{align}

Although \eqref{K} is a linear BSDE, its coefficients are unbounded since so is $\Lambda(i)$ (hence $\Gamma(i)$). And the equations in \eqref{K} are coupled through the term \break$``\sum_{j=1}^\ell q_{ij}K(j)"$.
Up to our knowledge, no existing literature could be directly applied to \eqref{K}. Next we will address ourselves to the solvability of \eqref{K} which is the main technique contribution of this paper.

\begin{remark}\label{oncoefficients1}
In our previous work \cite{HSX}, we studied a similar linear BSDE with unbounded coefficients, that is the BSDE for $(K,L)$ after Remark 5.9. By the connection $K = P H $ and $L = P \eta + K$, one can reduce the solvability issue of the BSDE for $(K,L)$ to that of
$(H,\eta)$ which satisfies \cite[(5.8)]{HSX}. Although (5.8) is still a linear BSDE with unbounded coefficients, the coefficient of $H$ becomes bounded so that it can be dealt by discounting.
The unbounded term can be removed by change of measure and the solvability issue can be resolved by contraction mapping. This is corresponding to the case that $a$ is bounded and $f=0$ in \eqref{Klinear} below.
\end{remark}

We first present several lemmas that will used to solve \eqref{K}.
The following lemma is called the John-Nirenberg inequality, which can be found in Theorem 2.2 of \cite{Ka}.
\begin{lemma}[John-Nirenberg Inequality]
\label{JohnNirenberg}
Suppose $\phi\in \BMO$ and \[\left\Vert\int_0^\cdot \phi(s)'dW(s)\right\Vert_{\mathrm{BMO}_2}<1.\]

Then for all $\{\mathcal{F}^W_t\}_{t\geq 0}$-stopping times $\tau\leq T$,
\begin{align*}
\E\Big[e^{\int_\tau^T|\phi(s)|^2ds}\;\Big|\;\mathcal{F}^W_\tau\Big]\leq \frac{1}{1-\left\Vert\int_0^\cdot \phi(s)'dW(s)\right\Vert_{\mathrm{BMO}_2}}.
\end{align*}
\end{lemma}
From this lemma, we immediately have the following estimate.
\begin{lemma}
\label{exponinteg}
Suppose $\phi\in \BMO$.
Then for any constants $p\in(0,2)$ and $K> 0$, there exists a constant $c_{p,K}>0$ such that
\begin{align*}
\E\Big[e^{K\int_\tau^T|\phi(s)|^pds}\;\Big|\;\mathcal{F}^W_\tau\Big]\leq c_{p,K}
\end{align*}
for all $\{\mathcal{F}^W_t\}_{t\geq 0}$-stopping times $\tau\leq T$.
\end{lemma}

\noindent \emph{Proof.} Denote $M(\cdot)=\int_0^\cdot \phi(s)'dW(s)$.
Let $\varepsilon$ be a constant such that $0<\varepsilon<\frac{1}{\left\Vert M\right\Vert_{\mathrm{BMO}_2}}$, then $\left\Vert\varepsilon M\right\Vert_{\mathrm{BMO}_2}=\varepsilon\left\Vert M\right\Vert_{\mathrm{BMO}_2}<1$.
For each $p\in(0,2)$ and $K>0$, we have $K|x|^p\leq \varepsilon^2x^2+c$, where
\[0<c=\sup_{x}\;(K|x|^p-\varepsilon^2x^2)<\infty.\]

Applying the John-Nirenberg inequality to $\varepsilon M$ yields
\begin{align*}
\E\Big[e^{K\int_\tau^T|\phi(s)|^p ds}\;\Big|\;\mathcal{F}^W_\tau\Big]&\leq \E\Big[e^{\int_\tau^T\big(\varepsilon^2|\phi(s)|^2+c\big)ds}\;\Big|\;\mathcal{F}^W_\tau\Big]\\
&\leq e^{cT} \E\Big[e^{\int_\tau^T\varepsilon^2|\phi(s)|^2ds}\;\Big|\;\mathcal{F}^W_\tau\Big]\\
&\leq \frac{e^{cT}}{1-||\varepsilon M||_{\mathrm{BMO}_2}}.
\end{align*}
This completes the proof.
\eof

The following lemma can be found in Page 26 of \cite{Ka}.
\begin{lemma}
\label{exponpower}
Suppose $\phi\in \BMO$. For any constant $p\geq1$, there is a generic constant $K_p>0$ such that
\begin{align*}
\E\left[\Big(\int_\tau^T|\phi(s)|^2ds\Big)^p\;\bigg|\;\mathcal{F}^W_{\tau}\right]\leq K_p
\left\Vert\int_0^\cdot \phi(s)'dW(s)\right\Vert_{\mathrm{BMO}_2}^{2p}
\end{align*}
for all $\{\mathcal{F}^W_t\}_{t\geq 0}$-stopping times $\tau\leq T$.
\end{lemma}

The following result solves a new class of one-dimensional BSDE with all coefficients being unbounded.
It will be used to establish the corresponding result in multi-dimensional case, that is, \eqref{K}.
\begin{lemma}
\label{onedimBSDE}
Suppose $p\in(0,2)$ is a constant and $\phi\in\BMO$. Suppose $a$ and $f$ are two $\mathbb{R}$-valued $\{\mathcal{F}^W_t\}_{t\geq 0}$-adapted processes, and $\beta$ is a $\mathbb{R}^n$-valued $\{\mathcal{F}^W_t\}_{t\geq 0}$-adapted process such that
\begin{align*}
|a|\leq |\phi|^p, \ |\beta|\leq|\phi|, \ |f|\leq |\phi|^2.
\end{align*}
Then for any $\xi\in L^\infty_{\mathcal{F}^W_T}(\Omega;\mathbb{R})$, the following $1$-dimensional BSDE
\begin{align}
\label{Y}
\begin{cases}
-dY=\big( a Y+ \beta' Z+ f\big)dt-Z'dW,\\
Y(T)=\xi,
\end{cases}
\end{align}
admits a unique solution $(Y,Z)\in L^{\infty}_{\mathcal{F}^W}(0,T;\mathbb{R})\times\BMO$.
\end{lemma}

\noindent \emph{Proof.} Introduce two processes
\begin{align*}
J(t)=\exp\left(\int_0^t a(s)ds\right),
\end{align*}
and
\begin{align*}
N(t)=\mathcal{E}\left(\int_0^t \beta(s)'dW(s)\right).
\end{align*}

Note that $N(t)$ is a uniformly integrable martingale, thus $\widetilde W(t):=W(t)-\int_0^t \beta(s)ds$ is a Brownian motion under the probability $\widetilde {\mathbb{P}}$ defined by
\begin{align*}
\frac{d\widetilde{ \mathbb{P}} }{d\mathbb{P}}\Bigg|_{\mathcal{F}^W_T}=N(T).
\end{align*}

Set $$Y(t)=J(t)^{-1}\widetilde{\E}_t\Big[J(T)\xi+ \int_t^TJ(s) f(s)ds\Big],$$ where $\widetilde\E$ is the expectation w.r.t. the probability measure $\widetilde{\mathbb{P}}$. Then clearly $Y(T)=\xi$.
Since $a$ is unbounded, so is $J$. Thus, we do not have the boundedness of $Y$ automatically. To show the boundedness of $Y$, we apply 
Lemmas \ref{exponinteg} and \ref{exponpower} to get
\begin{align*}
|Y(t)|&\leq \widetilde{\E}_t\Big[e^{\int_t^T a(r)dr}|\xi| + \int_t^T e^{\int_t^s a(r)dr}|f(s)|ds\Big]\\
&\leq\widetilde{\E}_t\Big[ce^{\int_t^T| \phi(r)|^pdr}\Big]+\widetilde{\E}_t\Big[e^{\int_t^T| \phi(r)|^pdr}\int_t^T| \phi(s)|^2ds\Big]\\
&\leq c+\frac{1}{2}\widetilde{\E}_t\Big[e^{2\int_t^T| \phi(r)|^pdr}+ \Big(\int_t^T| \phi(s)|^2ds\Big)^2\Big] \\
&\leq c,
\end{align*}
so $Y\in L^{\infty}_{\mathcal{F}^W}(0,T;\mathbb{R})$. Similarly,
\begin{align*}
&\quad\widetilde{\E}\Big(J(T)\xi+\int_0^TJ(s) f(s)ds\Big)^2\\
&\leq \widetilde{\E}\Big(J(T)|\xi|+\int_0^TJ(s)| f(s)|ds\Big)^2\\
&\leq 2c\widetilde{\E}\Big(e^{2\int_0^T| a(r)|dr}\Big)+2\widetilde{\E}\Big[\Big(\int_0^T| f(s)|e^{\int_0^T|a(r)|dr}ds\Big)^2\Big]\\
&\leq 2c\widetilde{\E}\Big(e^{2\int_0^T| \phi(r)|^pdr}\Big)+2\widetilde{\E}\Big[\Big(\int_0^T| \phi(s)|^2e^{\int_0^T|\phi(r)|^pdr}ds\Big)^2\Big]\\
&\leq c+2\widetilde{\E}\Big[e^{2\int_0^T|\phi(r)|^pdr}\Big(\int_0^T|\phi(s)|^2ds\Big)^2\Big]\\
&\leq c+ \widetilde{\E}\Big[e^{4\int_0^T| \phi(r)|^pdr}+\Big(\int_0^T| \phi(s)|^2ds\Big)^4\Big]\\
&<\infty.
\end{align*}

Thus
\[
J(t)Y(t)+\int_0^tJ(s)f(s)ds=\widetilde{\E}_t\Big[J(T)\xi+\int_0^TJ(s) f(s)ds\Big]
\]
is a square integrable martingale under $\widetilde{\mathbb{P}}$. By the martingale representation theorem, there exists $\widetilde Z\in L^2_{\mathcal{F}^{\widetilde W}}(0,T;\mathbb{R}^n)$ such that
\begin{align*}
J(t)Y(t)+\int_0^tJ(s) f(s)ds
&=J(0)Y(0)+\int_0^t\widetilde Z(s)'d\widetilde W(s).
\end{align*}

As a consequence,
\begin{align*}
d(J(t)Y(t)) &=\widetilde Z(t)'d\widetilde W(t)-J(t) f(t)dt\\
&= J(t)[Z(t)'d W(t)-Z(t)'\beta(t)dt- f(t)dt],
\end{align*}
where $Z(t)=J(t)^{-1}\widetilde Z(t)$.
By It\^o's lemma, 
\begin{align*}
dY(t)&= d(J(t)^{-1}\cdot J(t)Y(t))\\
&=(-a(t)J(t)^{-1})J(t)Y(t)dt+J(t)^{-1}(\widetilde Z(t)'d\widetilde W(t)-J(t)f(t)dt)\\
&=- a(t)Y(t)dt + Z(t)'d W(t)-Z(t)'\beta(t)dt- f(t)dt.
\end{align*}
Thus $(Y,Z)$ satisfies \eqref{Y}. Because $Y$ and $\xi$ are essentially bounded, using Lemma \ref{exponpower}, we get
\begin{align*}
\widetilde{\E}_{\tau}\Big[\int_{\tau}^T|Z(s)|^2ds\Big]&=\widetilde{\E}_{\tau}\Big[\Big(\int_{\tau}^TZ(s)'d\widetilde W(s)\Big)^2\Big]\\
&=\widetilde{\E}_{\tau}\Big[\Big(Y({\tau})-\xi-\int_{\tau}^T\big( a(s) Y(s)+ f(s)\big)ds\Big)^2\Big]\\
&\leq c+ c\widetilde{\E}_{\tau}\Big[\Big(\int_{\tau}^T ( |\phi(s)|^p + |\phi(s)|^2)ds\Big)^2\Big]\\
&\leq c+ c\widetilde{\E}_{\tau}\Big[\Big(\int_{\tau}^T ( 1 + 2|\phi(s)|^2)ds\Big)^2\Big]\\
&\leq c+ c\widetilde{\E}_{\tau}\Big[\Big(T+2\int_{\tau}^T |\phi(s)|^2ds\Big)^2\Big]\\
&\leq c,
\end{align*}
for all stopping times $\tau\leq T$. Hence, $\int_0^\cdot Z(s)'d\widetilde W(s)$ is a BMO martingale under $\widetilde {\mathbb{P}}$. Consequently $\int_0^\cdot Z(s)'d W(s)$ is a BMO martingale under ${\mathbb{P}}$. This shows that $(Y,Z)\in L^{\infty}_{\mathcal{F}^W}(0,T;\mathbb{R})\times\BMO$ is a solution of the $1$-dimensional BSDE \eqref{Y}.

Let us prove the uniqueness. Suppose
$$(Y,Z), \ (\hat Y, \hat Z) \in L^{\infty}_{\mathcal{F}^W}(0,T;\mathbb{R})\times\BMO$$
are both solutions of \eqref{Y}. Set
\[
\Delta Y=Y-\hat Y, \ \Delta Z=Z-\hat Z.
\]
Then $(\Delta Y, \Delta Z)$ satisfies the following BSDE:
\begin{align*}
\Delta Y(t)=\int_t^T\big( a(s) \Delta Y(s)+ \beta(s)' \Delta Z(s)\big)dt-\int_t^T\Delta Z(s)'dW(s).
\end{align*}
By It\^o's lemma, it follows
\begin{align*}
J(t)\Delta Y(t)=-\int_t^TJ(s)\Delta Z(s)'d\widetilde W(s).
\end{align*}

We get $\Delta Y=0$ by taking conditional expectation $\widetilde{\E}_t$ on both sides and using $J>0$. Thus
\begin{align*}
\widetilde{\E}_t\Big[\int_t^T(J(s)\Delta Z(s))^2ds\Big]&=\widetilde{\E}_t\Big[\Big(\int_t^TJ(s) \Delta Z(s)'d\widetilde W(s)\Big)^2\Big]\\
&=\widetilde{\E}_t\Big[(J(t)\Delta Y(t))^2\Big]=0,
\end{align*}
so $\Delta Z=0$ as $J>0$.
This completes the proof of the uniqueness.
\eof

%


With the help of the above 1-dimensional result and contraction mapping, we can solve a system of multi-dimensional BSDEs with all coefficients being unbounded.

\begin{theorem}\label{linearBSDEg}
Suppose $p\in(0,2)$ is a constant. Suppose, for every $i,j\in\cM$, $\phi(i)\in\BMO$, $f(i)$ and $\kappa_{ij}$ are $\mathbb{R}$-valued $\{\mathcal{F}^W_t\}_{t\geq 0}$-adapted processes, and $\beta(i)$ is a $\mathbb{R}^n$-valued $\{\mathcal{F}^W_t\}_{t\geq 0}$-adapted process such that
\begin{align*}
|\beta(i)|\leq|\phi(i)|, \ |f(i)|\leq |\phi(i)|^2, \ \sum_{j\in\cM}|\kappa_{ij}|\leq |\phi(i)|^p.
\end{align*}

Then for any given terminal value $(\xi(1),...,\xi(\ell))'\in L^\infty_{\mathcal{F}^W}(0,T;\mathbb{R}^{\ell})$,
the following system of (multi-dimensional) BSDEs:
\begin{align}\label{Klinear}
\begin{cases}
-dK(i)=\Big[\beta(i)' L(i)+ f(i)+\sum\limits_{j\in\cM} \kappa_{ij}K(j)\Big]dt-L(i)'dW,\\
K(T,i)=\xi(i), \ \mbox{for all} \ i\in\cM,
\end{cases}
\end{align}
admits a unique solution $(K(i),L(i))_{i\in\cM}$ such that
\[
(K(i),L(i))\in L^{\infty}_{\mathcal{F}^W}(0,T;\mathbb{R})\times\BMO, \ \mbox{for all} \ i\in\cM.
\]
\end{theorem}
\begin{remark}\label{oncoefficients2}
We emphasis here again that all of $\beta(i)$, $f(i)$ and $\kappa_{ij}$ are unbounded, whereas in \cite[(5.8)]{HSX} only $\beta(i)$ is unbounded.
\end{remark}

\noindent \emph{Proof.} For each $i\in\cM$, we introduce the process
\begin{align*}
N(t,i)=\mathcal{E}\left(\int_0^t\beta(s,i)'dW(s)\right).
\end{align*}

Note that $N(t,i)$ is a uniformly integrable martingale, thus $\widetilde W^i(t):=W(t)-\int_0^t\beta(s,i)ds$ is a Brownian motion under the probability $\widetilde {\mathbb{P}}^i$ defined by
\begin{align*}
\frac{d\widetilde{ \mathbb{P}}^i }{d\mathbb{P}}\Bigg|_{\mathcal{F}^W_T}=N(T,i), \ \mbox{ for all} \ i\in\cM.
\end{align*}

By Lemma \ref{onedimBSDE}, for any $U=(U(1),...,U(\ell))'\in L^\infty_{\mathcal{F}^W}(0,T;\mathbb{R}^{\ell})$ and each $i\in\cM$, the following $1$-dimensional linear BSDE admits a unique solution $(K(i),L(i))\in L^\infty_{\mathcal{F}^W}(0,T;\mathbb{R})\times\BMO$:
\begin{align*}
\begin{cases}
-dK(i)=\Big[\beta(i)' L(i)+f(i)+\sum\limits_{j\in\cM} \kappa_{ij}U(j)\Big]dt-L(i)'dW,\\
K(i,T)=\xi(i).
\end{cases}
\end{align*}

We let $\Theta$ denote the map $U\mapsto K:=(K(1),...,K(\ell))'$.
%

Thanks to Lemma \ref{exponpower}, there exists a constant $c_3>0$, independent of $t$ and $i$, such that
\begin{align}\label{c3}
&\quad\;\widetilde{\E}^i_t\bigg[ \bigg(\int_t^T |\phi(s,i)|^2ds\bigg)^{\frac{p}{2}} \bigg] <c_{3}.
\end{align}

For $U=(U(1),...,U(\ell))'\in L^\infty_{\mathcal{F}^W}(0,T;\mathbb{R}^{\ell})$, we introduce a new norm
\[
|U|_\infty:=\max\limits_{i\in\cM}\underset{(t, \omega)\in[0, T]\times\Omega}{\esssup}e^{c_4t}|U(t, i)|,
\]
where $c_{4}$ is a large positive constant to be determined.
Let $\mathcal{B}$ be the set of $U\in L^\infty_{\mathcal{F}^W}(0,T;\mathbb{R}^{\ell})$ with $|U|_\infty<\infty$.
For any $U$, $\widetilde U\in \mathcal{B}$, set $K=\Theta(U)$, $\widetilde K=\Theta(\widetilde U)$, and
\[\Delta K(t, i)=K(t, i)-\widetilde K(t, i), \ \text{and} \ \Delta U(t, i)=U(t, i)-\widetilde U(t, i).
\]

Then by It\^o's lemma,
\begin{align*}
\Delta K(t,i)=\widetilde{\E}^i_t\Big[\int_t^T\sum\limits_{j\in\cM}\kappa_{ij}(s)\Delta U(s,j)ds\Big].
\end{align*}

Since $p\in(0,2)$, it follows from H\"older's inequality that
\begin{align*}
e^{c_4t}|\Delta K(t,i)|&\leq e^{c_4t}\widetilde{\E}^i_t\Big[\int_t^T e^{-c_4s} \sum\limits_{j\in\cM}|\kappa_{ij}(s)|e^{c_4s}|\Delta U(s,j)|ds\Big]\\
&\leq \widetilde{\E}^i_t\Big[\int_t^T e^{-c_4(s-t)} |\phi(s,i)|^pds\Big]|\Delta U|_{\infty}\\
&\leq\widetilde{\E}^i_t\bigg[\bigg(\int_t^T e^{-\frac{2c_4}{2-p}(s-t)} ds\bigg)^{\frac{2-p}{2}}
\bigg(\int_t^T |\phi(s,i)|^2ds\bigg)^{\frac{p}{2}} \bigg]|\Delta U|_{\infty}\\
&\leq\widetilde{\E}^i_t\bigg[\bigg(\int_t^{\infty} e^{-\frac{2c_4}{2-p}(s-t)} ds\bigg)^{\frac{2-p}{2}}
\bigg(\int_t^T |\phi(s,i)|^2ds\bigg)^{\frac{p}{2}} \bigg]|\Delta U|_{\infty}\\
&= \Big(\frac{2-p}{2c_{4}}\Big)^{\frac{2-p}{2}}
\widetilde{\E}^i_t\bigg[ \bigg(\int_t^T |\phi(s,i)|^2ds\bigg)^{\frac{p}{2}} \bigg]|\Delta U|_{\infty}\\
&\leq \Big(\frac{2-p}{2c_{4}}\Big)^{\frac{2-p}{2}}c_3|\Delta U|_{\infty},
\end{align*}
where the last inequality is due to \eqref{c3}.
Let $c_{4}$ be sufficiently large such that $\Big(\frac{2-p}{2c_{4}}\Big)^{\frac{2-p}{2}}c_3\leq \frac{1}{2}$, then we have
\[
|\Delta K|_{\infty}
\leq \frac{1}{2}|\Delta U|_{\infty}.
\]

Therefore, $\Theta$ is a strict contraction mapping on $\mathcal{B}$ endowed with the norm $|\cdot|_\infty$.
Because $(\mathcal{B},\; |\cdot|_\infty)$ is a complete metric space, the map $\Theta$ admits a unique fixed point which is the unique solution to the $\ell$-dimensional BSDE \eqref{Klinear}.
\eof

\begin{corollary}
Under Assumptions \ref{assu1} and \ref{assu2}, the system of linear BSDEs \eqref{K} admits a unique solution $(K(i),L(i))_{i\in\cM}$ such that
\[
(K(i),L(i))\in L^{\infty}_{\mathcal{F}^W}(0,T;\mathbb{R})\times\BMO, \ \mbox{for all} \ i\in\cM.
\]
\end{corollary}

\noindent \emph{Proof.} Set
\begin{align*}
a(i)&= A(i)-B(i)'\Gamma(i),\\
\beta(i)&=C(i)-D(i)\Gamma(i),\\
f(i)&=q(i)Q(i)-P(i)(C(i)'\rho(i)+b(i))-\rho(i)'\Lambda(i)\\
&\quad\;+(P(i)D(i)'\rho(i)-R(i)p(i))'\Gamma(i),
\end{align*}
then $|a(i)|\leq c(1+|\Lambda(i)|), \ |\beta(i)|\leq c(1+|\Lambda(i)|), \ |f(i)|\leq c(1+|\Lambda(i)|)$, for all $i\in\cM$.

Also, \eqref{K} can be rewritten as
\begin{align*}
\begin{cases}
-dK(i)=\Big[a(i) K(i)+\beta(i)' L(i)+ f(i)+\sum\limits_{j\in\cM}q_{ij}K(j)\Big]dt-L(i)'dW,\\
K(T,i)=G(i)g(i), \ \mbox{for all} \ i\in\cM,
\end{cases}
\end{align*}
whence admits a unique solution $(K(i),L(i))_{i\in\cM}$ such that
\[
(K(i),L(i))\in L^{\infty}_{\mathcal{F}^W}(0,T;\mathbb{R})\times\BMO, \ \mbox{for all} \ i\in\cM,
\]
as a consequence of Theorem \ref{linearBSDEg} and $c(1+|\Lambda(i)|)\in\BMO$.
\eof
\section{Solution to the LQ problem \eqref{LQ}}\label{sec:solution}
The solution to the LQ problem \eqref{LQ} is stated as follows.
\begin{theorem}\label{LQcontrol}
Suppose that Assumptions \ref{assu1} and \ref{assu2} hold.
Let $(P(t,i), \ \Lambda(t,i))_{i\in\cM}$ and $(K(t,i), \ L(t,i))_{i\in\cM}$ be the unique solutions of the systems of BSDEs \eqref{P1} and \eqref{K}, respectively. Then the LQ problem \eqref{LQ} has an optimal control, as a feedback function of the time $t$, the state $X$, and the market regime $i$,
\begin{align}
\label{opticon}
u^*(t, X, i)&=-\big(R(t,i)+P(t,i)D(t,i)'D(t,i)\big)^{-1}\nonumber\\
&\quad\times
\Big[\Big(P(t,i)D(t,i)'C(t,i) +P(t,i)B(t,i)+D(t,i)'\Lambda(t,i)\Big)X\nonumber\\
&\quad+P(t,i)D(t,i)'\rho(t,i)-K(t,i)B(t,i)-D(t,i)'L(t,i)-R(t,i)p(t,i)\Big].
\end{align}

Moreover, the corresponding optimal value is
\begin{align}\label{optival}
V(x,i_0)&=P(0,i_0)x^2-2K(0,i_0)x+\E[G(\alpha_T)g(\alpha_T)^2]\nonumber\\
&\quad+\E\int_0^T\Big[P(t,\alpha_t)\rho(t,\alpha_t)'\rho(t,\alpha_t)-2K(t,\alpha_t)b(t,\alpha_t)-2\rho(t,\alpha_t)'L(t,\alpha_t)\nonumber\\
&\quad+Q(t,\alpha_t)q(t,\alpha_t)^2+p(t,\alpha_t)'R(t,\alpha_t)p(t,\alpha_t)\nonumber\\
&\quad-\big[D(t,\alpha_t)'(P(t,\alpha_t)\rho(t,\alpha_t)-L(t,\alpha_t))-K(t,\alpha_t)B(t,\alpha_t)\nonumber\\
&\quad-R(t,\alpha_t)p(t,\alpha_t)\big]'\big(R(t,\alpha_t)+P(t,\alpha_t)D(t,\alpha_t)'D(t,\alpha_t)\big)^{-1}\nonumber\\
&\quad\times\big[D(t,\alpha_t)'(P(t,\alpha_t)\rho(t,\alpha_t)-L(t,\alpha_t))-K(t,\alpha_t)B(t,\alpha_t)\nonumber\\
&\quad-R(t,\alpha_t)p(t,\alpha_t)\big]\Big]dt.
\end{align}
\end{theorem}

\noindent \emph{Proof.} The admissibility of the control process $u^*(t, X(t), \alpha_{t})$ will be proved in the following lemma.
The reminder of the proof is similar to that of Theorem 4.2 of \cite{HSX} via applying It\^{o}'s Lemma to $P(t,\alpha_{t})X(t)^2-2K(t,\alpha_{t})X(t)$, so we leave the details to the diligent readers.
\eof

\begin{lemma}
Under the conditions of Theorem \ref{LQcontrol}, we have $u^*(t, X(t), \alpha_{t})\in L^2_{\mathcal{F}}(0,T;\mathbb{R}^m)$.
\end{lemma}

\noindent \emph{Proof.} In light of the length of many equations, $``(t, X(t), \alpha_{t})"$ will be suppressed when no confusion occurs in the sequel.
Substituting \eqref{opticon} into the state process \eqref{state} (with $``i"$ replaced by $``\alpha_t"$), we have
\begin{align}
\label{SDEstate}
\begin{cases}
dX=\Big[\big[A-B'(R+PD'D)^{-1}(PD'C+PB+D'\Lambda)\big]X\\
\qquad\quad-B'(R+PD'D)^{-1}(PD'\rho-KB-D'L-Rp)+b\Big]dt\\
\qquad\quad+\Big[\big[C-D(R+PD'D)^{-1}(PD'C+PB+D'\Lambda)\big]X\\
\qquad\quad-D(R+PD'D)^{-1}(PD'\rho-KB-D'L-Rp)+\rho\Big]'dW\\
X(0)=x_0, \ \alpha_0=i_0.
\end{cases}
\end{align}

By the basic theorem on PP. 756-757 of Gal'chuk \cite{Ga}, the SDE \eqref{SDEstate} admits a unique strong solution.
For $(P(i), \Lambda(i))_{i\in\cM}$, $(K(i), L(i))_{i\in\cM}$, the unique solutions of \eqref{P1} and \eqref{K} respectively, and $X(t)$, the solution of \eqref{SDEstate},
applying It\^{o}'s lemma to $P(t, \alpha_t)X(t)^2-2K(t, \alpha_t)X(t)$, we have
\begin{align*}
& \quad\;\int_0^t\Big(Q(X-q)^2+(u^*-p)'R(u^*-p)\Big)dt+P(t,\alpha_t)X(t)^2-2K(t,\alpha_t)X(t)\\
&=P(0,i_0)x^2-2K(0,i_0)x+\int_0^t\Big[P\rho'\rho-2Kb-2L'\rho+Qq^2+p'Rp\\
& \quad -(PD'\rho-KB-D'L-Rp)'(R+PD'D)^{-1}(PD'\rho-KB-D'L-Rp)\Big]ds\\
& \quad +\int_0^t\Big[2(PX-K)(CX+Du^*+\rho)+X^2\Lambda-2XL\Big]'dW\\
&\quad+\int_0^t\Big\{X^2\sum_{j, j'\in\mathcal{M}}(P(s, j)-P(s, j'))I_{\{\alpha_{s-}=j'\}}\\
&\quad-2X\sum_{j, j'\in\mathcal{M}}(K(s, j)-K(s, j'))I_{\{\alpha_{s-}=j'\}}\Big\}d\tilde N_s^{j'j},
\end{align*}
where $(N^{j'j})_{j'j\in\mathcal{M}}$ are independent Poisson processes each with intensity $q_{j'j}$, and $\tilde N_t^{j'j}=N_t^{j'j}-q_{j'j}t, \ t\geq0$ are the corresponding compensated Poisson martingales under the filtration $\mathcal{F}$.

Because $X(t)$ is continuous, the stochastic integrals in the last equation are local martingales. Thus there exists an increasing sequence of stopping times $\tau_k$ such that $\tau_k\uparrow+\infty$ as $k\rightarrow+\infty$ such that
\begin{align}\label{L2con}
& \quad\; \E\Big[\int_0^{\iota\wedge\tau_k}\Big(Q(X-q)^2+(u^*-p)'R(u^*-p)\Big)ds\nonumber\\
&\qquad+P({\iota\wedge\tau_k})X({\iota\wedge\tau_k})^2-2K({\iota\wedge\tau_k})X({\iota\wedge\tau_k})\Big]\nonumber\\
&=P(0,i_0)x^2-2K(0,i_0)x+\E\int_0^{\iota\wedge\tau_k}\Big[P\rho'\rho-2Kb-2L'\rho+Qq^2+p'Rp\nonumber\\
& \quad-(PD'\rho-KB-D'L-Rp)'(R+PD'D)^{-1}(PD'\rho\nonumber\\
&\quad-KB-D'L-Rp)\Big]ds
\end{align}
for all stopping times $\iota\leq T$.

Under Assumptions \ref{assu1} and \ref{assu2} (i), we have\newpage
\begin{align*}
& \quad\; \E\Big[\int_0^{\iota\wedge\tau_k}\delta|u^*-p|^2ds\Big]\\
&\leq \E\Big[\int_0^{\iota\wedge\tau_k}(u^*-p)'R(u^*-p)ds
+P({\iota\wedge\tau_k})\Big(X({\iota\wedge\tau_k})-\frac{K(t,\alpha_t)}{P(t,\alpha_t)}\Big)^2\Big]\\
&\leq P(0,i_0)x^2-2K(0,i_0)x+\E\int_0^{\iota\wedge\tau_k}\Big(P\rho'\rho-2Kb-2L'\rho+Qq^2+p'Rp\Big)ds\\
&\quad+\E\Big[\frac{K(\iota\wedge\tau_k)^2}{P(\iota\wedge\tau_k)}\Big]\\
&\leq P(0,i_0)x^2-2K(0,i_0)x+\E\int_0^{T}\Big|P\rho'\rho-2Kb-2L'\rho+Qq^2+p'Rp\Big|ds\\
&\quad+\E\Big[\frac{K(\iota\wedge\tau_k)^2}{P(\iota\wedge\tau_k)}\Big]\\
&\leq c,
\end{align*}
where the constant $c>0$ is independent of $k$. Taking $\iota=T$ and letting $k\rightarrow\infty$, it follows from the monotone theorem and the boundedness of $p$ that
$$u^*(t, X(t), \alpha_{t})\in L^2_{\mathcal{F}}(0,T;\mathbb{R}^m).$$

Similarly, under Assumptions \ref{assu1} and \ref{assu2} (ii), we have
\begin{align*}
& \quad\; \E\Big[P({\iota\wedge\tau_k})\Big(X({\iota\wedge\tau_k})-\frac{K(\iota\wedge\tau_k)}{P(\iota\wedge\tau_k)}\Big)^2\Big]\\
&\leq P(0,i_0)x^2-2K(0,i_0)x+\E\int_0^{T}\Big|P\rho'\rho-2Kb-2L'\rho+Qq^2+p'Rp\Big|ds\\
&\quad+\E\Big[\frac{K(\iota\wedge\tau_k)^2}{P(\iota\wedge\tau_k)}\Big]\\
&\leq c.
\end{align*}

By Lemma \ref{Riccatisingular}, there exists a constant $c_4>0$ such that $P(i)\geq c_4$, for all $i\in\cM$. Therefore
\begin{align*}
c_4\E\Big[X({\iota\wedge\tau_k})^2\Big]&\leq \E\Big[P({\iota\wedge\tau_k})X({\iota\wedge\tau_k})^2\Big]\\
&\leq 2\E\Big[P({\iota\wedge\tau_k})\Big(X({\iota\wedge\tau_k})-\frac{K(\iota\wedge\tau_k)}{P(\iota\wedge\tau_k)}\Big)^2\Big]
+2\E\Big[\frac{K(\iota\wedge\tau_k)^2}{P(\iota\wedge\tau_k)}\Big]\leq c.
\end{align*}

Letting $k\rightarrow\infty$, it follows from Fatou's lemma that
\begin{align*}
\E\Big[X(\iota)^2\Big]\leq c,
\end{align*}
for all stopping times $\iota\leq T$.
This further implies
\begin{align}\label{bound2}
\E \int_0^{\iota\wedge T}X(s)^2ds\leq \int_0^{T}\E\left[X(s)^2\right]ds\leq cT.
\end{align}

By It\^{o}'s Lemma, we have
\begin{align*}
X(t)^2&=x^2+\int_0^t\Big[(u^*)'D'Du^*+2X(D'C+B)'u^*+2\rho'Du^*\\
&\quad+(2A+C'C)X^2+2X(b+C'\rho)+\rho'\rho\Big]ds\\
&\quad\;+\int_0^t2X(CX+Du^*+\rho)'dW.
\end{align*}

Assumptions \ref{assu1} and \ref{assu2} (ii) and the positiveness of $P$ implies $R+PD'D\geq c>0$, so
\begin{align}\label{uestimate}
|u^*|\leq c(1+ |\Lambda||X|+|L|).
\end{align}

Let
\[\theta_k=\inf\Big\{t\geq 0: |X(t)|+\int_0^t (|\Lambda(s)|^2+|L(s)|^2)ds>k\Big\}. \]

Because $X(t)$ is continuous, $\Lambda$, $L\in L^{2}_{\mathcal {F}^W}(0, T;\mathbb{R}^n)$, it follows that
\begin{align*}
&\ \ \ \ x^2+\E \int_0^{T\wedge\theta_k}(u^*)'D'Du^*ds\\
&=\E\Big[ X(T\wedge\theta_k)^2\Big]-\E \int_0^{T\wedge\theta_k} \Big[2X(D'C+B)'u^*+2\rho'Du^*\\
&\quad+(2A+C'C)X^2+2X(b+C'\rho)+\rho'\rho\Big]ds.
\end{align*}

Let $\delta>0$ be given in Assumption \ref{assu2}. By Assumption \ref{assu1} and \eqref{bound2}, the above by the elementary inequality $2ab\leq \frac{\epsilon}{2}a^2+\frac{2}{\epsilon}b^2$ leads to
\begin{align*}
&\quad\;\delta \E \int_0^{T\wedge\theta_k}|u^*(s, X(s), \alpha_s)|^2ds\\
&\leq c+c\E \int_0^{T\wedge\theta_k}\Big[X(s)^2+2(|X(s)|+1) |u^*|\Big]ds\\
&\leq c+\frac{4c^2}{\delta}+(c+\frac{4c^2}{\delta})\E \int_0^{T\wedge\theta_k}X(s)^2ds+\frac{\delta}{2}\E \int_0^{T\wedge\theta_k} |u^*(s, X(s), \alpha_s)|^{2}ds\\
&\leq c+\frac{\delta}{2}\E \int_0^{T\wedge\theta_k} |u^*(s, X(s), \alpha_s)|^{2}ds.
\end{align*}

The last expectation can be shown to be finite by \eqref{uestimate} and the definition of $\theta_k$.
So, after rearrangement and sending $k\to\infty$, it follows from the monotone convergence theorem that
$$u^*(t, X(t), \alpha_{t})\in L^2_{\mathcal{F}}(0,T;\mathbb{R}^m).$$

The proof is complete. \eof

%
%

\begin{remark}
Set $(h(t,i),\eta(t,i))=\Big(\frac{K(t,i)}{P(t,i)},
-\frac{K(t,i)\Lambda(t,i)}{P(t,i)^2}+\frac{L(t,i)}{P(t,i)}\Big)$, then\break $(h(t,i),\eta(t,i))_{i\in\cM}$ is the solution to the following system of ($\ell$-dimensional) linear BSDEs
\begin{align}
\label{h}
\begin{cases}
dh(i)=\Big\{\Big[A(i)+C(i)'C(i)+C(i)'\frac{\Lambda(i)}{P(i)}+\frac{Q(i)}{P(i)}-
C(i)'D(i)\Gamma(i)\Big]h(i)\\
\qquad\qquad-\Big[C(i)+\frac{\Lambda(i)}{P(i)}-
D(i)\Gamma(i)\Big]'\eta(i)-(D(i)'\rho(i)-\frac{R(i)p(i)}{P(i)})'\Gamma(i)\\
\qquad\qquad-\frac{q(i)Q(i)}{P(i)}+b(i)+\rho(i)'C(i)+\frac{\rho(i)'\Lambda(i)}{P(i)}\\
\qquad\qquad+\frac{1}{P(i)}\sum\limits_{j=1}^\ell q_{ij}P(j)(h(i)-h(j))\Big\}dt+\eta(i)'dW,\\
h(T,i)=g(i), \ \mbox{for all} \ i\in\cM.
\end{cases}
\end{align}

Applying It\^{o}'s Lemma to $K(t,\alpha_{t})h(t,\alpha_t)$ on $[0,T]$, the optimal value \eqref{optival} could be represented by $(h(t,i),\eta(t,i))_{i\in\cM}$:
\begin{align*}
V(x,i_0)&=P(0,i_0)(x-h(0,i_0))^2+\E\int_0^T Q(h-q)^2dt\\
&\quad+\E\int_0^TP\Big(\rho+hC-\eta\Big)'
\Big(I_{n}-PD(R+PD'D)^{-1}D'\Big)\Big(\rho+hC-\eta\Big)dt\\
&\quad+\E\int_0^T\Big[p'(R-R(R+PD'D)^{-1}R)p\\
&\quad+2P\Big(\rho+hC-\eta\Big)'D(R+PD'D)^{-1}Rp\Big]dt\\
&\quad+\E\int_0^T\sum_{j=1}^\ell q_{\alpha_t j}P(t,j)\Big(h(t,\alpha_t)-h(t,j)\Big)^2dt\\
&=P(0,i_0)(x-h(0,i_0))^2\\
&\quad+\E\int_0^T\Big[ Q(h-q)^2dt+P|\rho+hC-\eta|^2\Big]
\Big(\rho+hC-\eta\Big)dt\\
&\quad+\E\int_0^T\Big[p'Rp-\Big(Rp-PD'(\rho+hc-\eta)\Big)'(R+PD'D)^{-1}\\
&\quad\times\Big(Rp-PD'(\rho+hc-\eta)\Big)\Big]dt\\
&\quad+\E\int_0^T\sum_{j=1}^\ell q_{\alpha_t j}P(t,j)\Big(h(t,\alpha_t)-h(t,j)\Big)^2dt,
\end{align*}
where $``(t, \alpha_{t})"$ are suppressed for simplicity.
\end{remark}


\section{Application to a mean-variance asset-liability management problem}\label{sec:application}
Consider a financial market consisting of a risk-free asset (the money market
instrument or bond) whose price is $S_{0}$ and $m$ risky securities (the
stocks) whose prices are $S_{1}, \ldots, S_{m}$. Assume $m\leq n$, i.e., the number of risky securities is no more than the dimension of the Brownian motion. The financial market is incomplete if $m<n$.
These asset prices are driven by stochastic differential equations (SDEs):
\begin{align*}
\begin{cases}
dS_0(t)=r(t,\alpha_t)S_{0}(t)dt, \\
S_0(0)=s_0>0,
\end{cases}
\end{align*}
and
\begin{align*}
\begin{cases}
dS_k(t)=S_k(t)\Big((\mu_k(t, \alpha_t)+r(t,\alpha_t))dt+\sum\limits_{j=1}^n\sigma_{kj}(t, \alpha_t)dW_j(t)\Big), \\
S_k(0)=s_k>0,
\end{cases}
\end{align*}
where $r(t,i)$ is the interest rate process and $\mu_k(t, i)$ and $\sigma_k(t, i):=(\sigma_{k1}(t, i), \ldots, $ $\sigma_{kn}(t, i))$ are the mean excess return rate process and volatility rate process of the $k$th risky security corresponding to a market regime $\alpha_t=i$, for every $k=1, \ldots, m$ and $i\in\cM$.

Define the mean excess return vector
\begin{align*}
\mu(t, i)=(\mu_1(t, i), \ldots, \mu_m(t, i))',
\end{align*}
and volatility matrix
\begin{align*}
\sigma(t, i)=
\left(
\begin{array}{c}
\sigma_1(t, i)\\
\vdots\\
\sigma_m(t, i)\\
\end{array}
\right)
\equiv (\sigma_{kj}(t, i))_{m\times n}, \ \text{for}\ \text{each} \ i\in\cM.
\end{align*}

A small investor, whose actions cannot affect the asset prices, needs to decide at every time
$t\in[0, T]$ the amount $\pi_j(t)$ to invest in the $j$th risky asset, $j=1, \ldots, m$. The vector process $\pi(\cdot):=(\pi_1(\cdot), \ldots, \pi_m(\cdot))'$ is called a portfolio of the investor.
The admissible portfolio set is defined as
\begin{align*}
\mathcal U=L^2_{\mathcal F}(0, T;\mathbb R^m).
\end{align*}

Then the investor's asset value $\ass(\cdot)$ corresponding to a portfolio $\pi(\cdot)$ is the unique strong solution of the SDE:
\begin{align}
\label{asset}
\begin{cases}
d\ass(t)=[r(t,\alpha_t)\ass(t)+\pi(t)'\mu(t, \alpha_t)]dt+\pi(t)'\sigma(t, \alpha_t)dW(t), \\
\ass(0)=\ass_0, \ \alpha_0=i_0.
\end{cases}
\end{align}

Besides the asset value above, the investor has to pay for some liability $l(\cdot)$ whose value is modeled as an It\^{o} process
\begin{align*}
\begin{cases}
dl(t)=[r(t,\alpha_t)l(t)-b(t,\alpha_t)]dt-\rho(t,\alpha_t)'dW(t),\\
l(0)=l_0, \ \alpha_0=i_0.
\end{cases}
\end{align*}

Then the surplus value of the investor $X(t):=\ass(t)-l(t)$ is governed by
\begin{align}
\label{wealth}
\begin{cases}
dX(t)=[r(t,\alpha_t)X(t)+\pi(t)'\mu(t, \alpha_t)+b(t,\alpha_{t})]dt\\
\qquad\qquad\qquad+[\pi(t)'\sigma(t, \alpha_t)+\rho(t,\alpha_{t})']dW(t), \\
X(0)=x:=\ass_0-l_0, \ \alpha_0=i_0.
\end{cases}
\end{align}

For a given expectation level $z\in\mathbb{R}$, the investor's mean-variance asset-liability management problem is
\begin{align}
\mathrm{Minimize}&\quad \mathrm{Var}(X(T))=\E\big[(X(T)-z)^{2}\big]%
, \nonumber\\
\mathrm{ s.t.} &\quad
\begin{cases}
\E (X(T))=z, \\
\pi\in \mathcal{U}.
\end{cases}
\label{optm}%
\end{align}
\par


\begin{remark}
The liability process is modeled as a geometric Brownian motion in \cite{CYY}, \cite{CL}, and as a Brownian motion with drift in \cite{Xie}, \cite{XLW},
\cite{ZL}.
As explained in \cite{ZL},
``The liability here is in a generalized sense. We understand it as the subtraction of the real
liability and the stochastic income of the investor... A negative liability means that the stochastic income of the investor is bigger than his/her real liability". We can also interpret the liability as the total value of the investor's non-tradable assets.
\end{remark}

We impose the following assumption.
\begin{assumption} \label{assu5}
\emph{For all $i\in\cM$,}
\begin{align*}
\begin{cases}
r(\cdot, \cdot, i), \ b(\cdot, \cdot, i)\in L_{\mathcal{F}^W}^\infty(0, T;\mathbb{R}), \ \mu(\cdot, \cdot, i)\in L_{\mathcal{F}^W}^\infty(0, T;\mathbb{R}^m), \\
\rho(\cdot, \cdot, i)\in L_{\mathcal{F}^W}^\infty(0, T;\mathbb{R}^n), \
\sigma(\cdot, \cdot, i)\in L_{\mathcal{F}^W}^\infty(0, T;\mathbb{R}^{m\times n}),
\end{cases}
\end{align*}
\emph{and $\sigma(t,i)\sigma(t,i)'\geq\delta I_{m}$ with some constant $\delta>0$, for a.e. $t\in[0,T]$.}
\end{assumption}

We shall say that the problem \eqref{optm} is feasible for a given $z$ if there is a portfolio $\pi\in \mathcal{U}$ which satisfies the target constraint $\E(X(T))=z$.
The following result gives necessary and sufficient conditions for the feasibility of (\ref{optm}) for any $z\in\mathbb{R}$.
\begin{theorem}
Suppose that Assumption \ref{assu5} holds.
Let $(\psi(t,i),\xi(t,i))\break\in $ $L^\infty_{\mathcal{F}^W}(0,T;\mathbb{R})\times L^2_{\mathcal{F}^W}(0,T;\mathbb{R}^n), \ i=1,...,\ell$ be the unique solution of system of linear BSDEs:
\begin{align}
\begin{cases}
d\psi(t,i)=-\Big(r(t,i)\psi(t,i)+\sum\limits_{j=1}^\ell q_{ij}\psi(t,j)\Big)dt+\xi(t,i)'dW_t,\\
\psi(T,i)=1, \ \mbox{ for all $i\in\cM$}.
\end{cases}
\end{align}

Then the mean-variance asset-liability management problem (\ref{optm}) is feasible for any $z\in\mathbb{R}$ if and only if
\begin{equation}
\label{feasible}
\E \int_0^T|\psi(t,\alpha_t)\mu(t,\alpha_t)+\sigma(t,\alpha_t)\xi(t,\alpha_t)|^2dt>0.
\end{equation}
\end{theorem}

\noindent \emph{Proof.} For any $\pi\in\mathcal{U}$ and any real number $\beta$, set a portfolio $\pi^{\beta}(t):=\beta\pi(t)$. Let $X^\beta$ be the wealth process corresponding to $\pi^\beta$. Then $X^\beta(t)=X^0(t)+\beta Y(t)$, where
\begin{align*}
\begin{cases}
dX^0(t)=[r(t,\alpha_t)X^0(t)+b(t,\alpha_t)]dt+\rho(t,\alpha_{t})'dW(t),\\
X^0(0)=x, \ \alpha_0=i_0,
\end{cases}
\end{align*}
and
\begin{align*}
\begin{cases}
dY(t)=[r(t,\alpha_t)Y(t)+\pi(t)'\mu(t,\alpha_t)]dt+\pi(t)'\sigma(t,\alpha_t)dW(t),\\
Y(0)=0, \ \alpha_0=i_0.
\end{cases}
\end{align*}

We first prove the ``if'' part.
Let $\pi(t)=\psi(t,\alpha_t)\mu(t,\alpha_t)+\sigma(t,\alpha_t)\xi(t,\alpha_t)$, then $\pi\in\mathcal{U}$.
Applying It\^o's lemma to $Y(t)\psi(t,\alpha_t)$, we have
\begin{align*}
\E (X(T))&=\E (X^0(T))+\beta \E (Y(T))\\
&=\E (X^0(T))+\beta\E \int_0^T\pi(t)'(\psi(t,\alpha_t)\mu(t,\alpha_t)+\sigma(t,\alpha_t)\xi(t,\alpha_t))dt\\
&=\E (X^0(T))+\beta\E \int_0^T|\psi(t,\alpha_t)\mu(t,\alpha_t)+\sigma(t,\alpha_t)\xi(t,\alpha_t)|^2dt.
\end{align*}

Notice that $\E(X^0(T))$ is a constant independent of $\pi$, then under (\ref{feasible}), for any $z\in\mathbb{R}$, there exists $\beta\in\mathbb{R}$ such that $\E (X(T))=z$.

Conversely, suppose that (\ref{optm}) is feasible for any $z\in\mathbb{R}$. Then for any $z\in\mathbb{R}$, there is a $\pi\in\mathcal U$, such that $\E (X(T))=\E (X^0(T))+ \E (Y(T))=z$. Notice that $\E (X^0(T))$ is independent of $\pi$, thus it is necessary that there is a $\pi\in\mathcal U$ such that $\E (Y(T))\neq0$.
It follows from $$\E (Y(T))=\E \int_0^T\pi(t)'(\psi(t,\alpha_t)\mu(t,\alpha_t)+\sigma(t,\alpha_t)\xi(t,\alpha_t))dt$$ that (\ref{feasible}) is true.
\eof

If \eqref{feasible} does not hold, the above proof shows that there is only one feasible target $z$. To avoid this trivial case, we assume \eqref{feasible} holds from now on. This allows us to deal with the constraint $\E(X(T))=z$ by Lagrangian method.

We introduce a Lagrange multiplier $-2\lambda\in\mathbb{R}$ and consider the following \emph{relaxed}
optimization problem:
\begin{align}\label{optmun}
\mathrm{Minimize} &\quad {\mathbb{E}}(X(T)-z)^{2}-2\lambda({\mathbb{E}}%
X(T)-z)={\mathbb{E}}(X(T)-(\lambda+z))^{2}-\lambda^{2}=:\hat{J}(\pi, \lambda), \\
\mathrm{s.t.} &\quad \pi\in \mathcal{U}.\nonumber%
\end{align}

Problems \eqref{optm} and \eqref{optmun} are linked by the Lagrange duality theorem (see Luenberger \cite{Lu})
\begin{align}\label{duality}
\min_{\pi\in\mathcal{U}, \E(X(T))=z}\mathrm{Var}(X(T)%
)=\max_{\lambda\in\mathbb{R}}\min_{\pi\in\mathcal{U}}\hat{J}(\pi, \lambda).
\end{align}

This allows us to solve the problem \eqref{optm} by a two-step procedure: First solve the relaxed problem \eqref{optmun}, then find a $\lambda^{*}$ to maximize $\lambda\mapsto\min_{\pi\in\mathcal{U}}\hat{J}(\pi, \lambda)$.



Apparently, the problem \eqref{optmun} is a special case of the LQ problem \eqref{LQ} where $A=r, \ B=\mu, \ C=0, \ D'=\sigma, \ Q=0, \ q=0, \ R=0, \ p=0, \ G=1$.
Recall that $``(t, \alpha_{t})"$ or $``(t, i)"$ are often suppressed where no confusion occurs for simplicity.
Furthermore, the system of BSDEs \eqref{P1} becomes
\begin{align}
\label{P2}
\begin{cases}
dP(i)=-\Big[2rP(i)-\frac{1}{P(i)}(P(i)\mu+\sigma\Lambda(i))'(\sigma\sigma')^{-1}(P(i)\mu+\sigma\Lambda(i))\\
\quad\quad\qquad+\sum\limits_{j=1}^\ell q_{ij}P(j)\Big]dt+\Lambda(i)'dW, \\
P(T,i)=1,\\
P(i,t)>0, \ \mbox{ for all $i\in\cM$},
\end{cases}
\end{align}
and \eqref{K} becomes
\begin{align}\label{K2}
\begin{cases}
d K(i)=-\Big[\big(r-\frac{1}{P(i)}(P(i)\mu+\sigma\Lambda(i))'(\sigma\sigma')^{-1}\mu\big) K(i)
\\
\qquad\qquad-\frac{1}{P(i)}(P(i)\mu+\sigma\Lambda(i))'(\sigma\sigma')^{-1}\sigma L(i)+(P(i)\mu+\sigma\Lambda(i))'(\sigma\sigma')^{-1}\sigma\rho\\
\qquad\qquad-P(i)b-\rho'\Lambda(i)+\sum\limits_{j=1}^\ell q_{ij} K(j)\Big]dt+ L(i)'dW,\\
K(T,i)=\lambda+z,\ \mbox{for all} \ i\in\cM.
\end{cases}
\end{align}

From Theorem \ref{LQcontrol}, we immediately have
\begin{theorem}\label{optalm}
Suppose that Assumption \ref{assu5} holds.
Let $(P(t,i), \ \Lambda(t,i))_{i\in\cM}$ and $(K(t,i), L(t,i))_{i\in\cM}$ be the unique solutions of \eqref{P2} and \eqref{K2}, respectively. Then the problem \eqref{optmun} has an optimal control, as a feedback function of the time $t$, the state $X$, and the market regime $i$,
\begin{align*}
\pi^*(t, X, i)&=-\frac{1}{P(t,i)}\big(\sigma(t,i)\sigma(t,i)'\big)^{-1}\Big[\Big(P(t,i)\mu(t,i)+\sigma(t,i)\Lambda(t,i)\Big)X\nonumber\\
&\quad+P(t,i)\sigma(t,i)\rho(t,i)-( K(t,i)\mu(t,i)+\sigma(t,i) L(t,i))\Big].
\end{align*}

Moreover, the corresponding optimal value is
\begin{align*}
\min_{\pi\in\mathcal{U}}\hat{J}(\pi, \lambda)&=P(0,i_0)x^2-2 K(0,i_0)x\\& \quad+(\lambda+z)^2-\lambda^2+\E\int_0^T\Big[P\rho'\rho-2( Kb+\rho' L)\\
&\quad-\frac{1}{P}(P\sigma\rho-( K\mu+\sigma L))'(\sigma\sigma')^{-1}(P\sigma\rho-( K\mu+\sigma L))\Big]dt.
\end{align*}

\begin{remark}
In this case, \eqref{h} becomes
\begin{align}
\label{h2}
\begin{cases}
d h(i)=\Big\{r h(i)+\mu'(\sigma\sigma')^{-1}\sigma\eta(i)+\frac{1}{P(i)}\Lambda(i)'(\sigma'(\sigma\sigma')^{-1}\sigma-I_n)\eta(i)\\
\qquad\qquad-\Big[\frac{1}{P(i)}(P(i)\mu+\sigma\Lambda(i))'(\sigma\sigma')^{-1}\sigma\rho-b-\frac{1}{P(i)}\rho'\Lambda(i)\Big]\\
\qquad\qquad
+\frac{1}{P(i)}\sum\limits_{j=1}^\ell q_{ij}P(j)( h(i)- h(j))\Big\}dt+\eta(i)'dW,\\
h(T,i)=\lambda+z, \ \mbox{for all} \ i\in\cM.
\end{cases}
\end{align}

The optimal control and optimal value in Theorem \ref{optalm} can be rewritten in terms of the unique solution $(h(t,i), \eta(t,i))_{i\in\cM}$ of \eqref{h2}:
\begin{align*}
\pi^*(t, X, i)=-\frac{1}{P}(\sigma\sigma')^{-1}\Big[(P\mu+\sigma\Lambda)(X- h)+P\sigma\rho-P\sigma \eta\Big],
\end{align*}
and
\begin{align}\label{optivalalm}
\min_{\pi\in\mathcal{U}}\hat{J}(\pi, \lambda)
&=P(0,i_0)(x- h(0,i_0))^2-\lambda^2\\&\quad+\E\int_0^T\sum_{j\in\cM} q_{\alpha_tj}P(j)( h(t,\alpha_t)- h(t,j))^2dt\nonumber\\
&\quad+\E\int_0^T\Big[P(\rho-\eta)'(I_{n}-\sigma'(\sigma\sigma')^{-1}\sigma)(\rho-\eta)\Big]dt.
\end{align}
\end{remark}
\end{theorem}

\begin{remark}
If $m=n=1$ and $r(\cdot, i), \ \mu(\cdot, i), \ \sigma(\cdot, i), \ b(\cdot, i), \ \rho(\cdot, i)$ are deterministic functions of $t$ for all $i\in\cM$, then $\Lambda(i)=L(i)=\eta(i)\equiv0$. Accordingly,
\eqref{P2} \eqref{K2}, \eqref{h2} and \eqref{optivalalm} degenerate to the ODEs (17), (18), (23) and Eq. (26) in Xie \cite{Xie} respectively with $\rho(t)\equiv 1$ (Here we take the notation of $\rho(t)$ used in \cite{Xie}).
\end{remark}

By the Lagrange duality relationship \eqref{duality}, we need to find $\lambda^*\in\mathbb{R}$ which attains the optimal value $\min_{\pi\in\mathcal{U}}\hat{J}(\pi, \lambda)$. Notice that $( h(t,i),\eta(t,i))_{i\in\cM}$ depends on $\lambda$, we need to seperate $\lambda$ from the equations of $( h(t,i),\eta(t,i))_{i\in\cM}$.

Let $( h_1(t,i),\eta_1(t,i))_{i\in\cM}$ and $( h_2(t,i),\eta_2(t,i))_{i\in\cM}$ be, respectively, the unique solutions of the following two systems of linear BSDEs,
\begin{align*}
\begin{cases}
d h_1(i)=\Big\{r h_1(i)+\mu'(\sigma\sigma')^{-1}\sigma\eta_1(i)+\frac{1}{P(i)}\Lambda(i)'(\sigma'(\sigma\sigma')^{-1}\sigma-I_n)\eta_1(i)\\
\qquad\qquad-\Big[\frac{1}{P(i)}(P(i)\mu+\sigma\Lambda(i))'(\sigma\sigma')^{-1}\sigma\rho-b-\frac{1}{P(i)}\rho'\Lambda(i)\Big]\\
\qquad\qquad
+\frac{1}{P(i)}\sum\limits_{j\in\cM}q_{ij}P(j)( h_1(i)- h_1(j))\Big\}dt+\eta_1(i)'dW,\\
h_1(T,i)=0, \ \mbox{for all} \ i\in\cM,
\end{cases}
\end{align*}
and
\begin{align*}
\begin{cases}
d h_2(i)=\Big\{r h_2(i)+\mu'(\sigma\sigma')^{-1}\sigma\eta_2(i)+\frac{1}{P(i)}\Lambda(i)'(\sigma'(\sigma\sigma')^{-1}\sigma-I_n)\eta_2(i)\\
\qquad\qquad
+\frac{1}{P(i)}\sum\limits_{j\in\cM}q_{ij}P(j)( h_2(i)- h_2(j))\Big\}dt+\eta_2(i)'dW,\\
h_2(T,i)=1, \ \mbox{for all} \ i\in\cM.
\end{cases}
\end{align*}

Then by uniqueness of the solution of \eqref{h2}, it is not hard to verify $h=h_1+(\lambda+z)h_2$ and $\eta=\eta_1+(\lambda+z)\eta_2$. For notation simplicity, we
denote
\begin{align*}
P_0:=P(0,i_0), \ h_{1,0}:=h_1(0,i_0), \ h_{2,0}:=h_2(0,i_0).
\end{align*}

Then from \eqref{optivalalm}, we have
\begin{align*}
\min_{\pi\in\mathcal{U}}\hat{J}(\pi, \lambda)
&=P_0\Big(x- h_{1,0}-(\lambda+z)h_{2,0}\Big)^2-\lambda^2\\
&\quad+\E\int_0^T\sum_{j\in\cM} q_{\alpha_tj}P(j)( h_1(\alpha_t)- h_1(j))^2dt\\
&\quad+(\lambda+z)^2\E\int_0^T\sum_{j\in\cM} q_{\alpha_tj}P(j)( h_2(\alpha_t)- h_2(j))^2dt\\
&\quad+2(\lambda+z)\E\int_0^T\sum_{j\in\cM} q_{\alpha_tj}P(j)( h_1(\alpha_t)- h_1(j))(h_2(\alpha_t)-h_2(j))dt\\
&\quad+\E\int_0^T\Big[P(\rho-\eta_1-(\lambda+z)\eta_2)'(I_{ n}-\sigma'(\sigma\sigma')^{-1}\sigma)\\
&\quad\times(\rho-\eta_1-(\lambda+z)\eta_2)\Big]dt\\
&=-(1-P_0h_{2,0}^2-M_1)\lambda^2+(P_{0}h_{2,0}^2+M_1)z^2\\
&\quad+2\Big(M_2+(P_{0}h_{2,0}^2+M_1)z-P_{0}h_{2,0}(x-h_{1,0})\Big)\lambda\\
&\quad+2(M_2-P_{0}h_{2,0})(x-h_{1,0}))z+M_3+P_0(x-h_{1,0})^2,
\end{align*}

\noindent where
\begin{align*}
M_1:&=\E\int_0^T\sum_{j\in\cM} q_{\alpha_tj}P(j)( h_2(\alpha_t)- h_2(j))^2dt
\\&\quad+\E\int_0^TP\eta_2'(I_{n}-\sigma'(\sigma\sigma')^{-1}\sigma)\eta_2dt,\\
M_2:&=\E\int_0^T\sum_{j\in\cM} q_{\alpha_tj}P(j)( h_1(\alpha_t)- h_1(j))(h_2(i)-h_2(j))dt\\
&\quad-\E\int_0^TP(\rho-\eta_1)'(I_{ n}-\sigma'(\sigma\sigma')^{-1}\sigma)\eta_2dt,\\
M_3:&=\E\int_0^T\sum_{j\in\cM} q_{\alpha_tj}P(j)( h_1(\alpha_t)- h_1(j))^2dt\\
&\quad+\E\int_0^TP(\rho-\eta_1)'(I_{ n}-\sigma'(\sigma\sigma')^{-1}\sigma)(\rho-\eta_1)dt.
\end{align*}

By Theorem 5.11 of \cite{HSX}, $0<P_0h_{2,0}^2+M_1<1$. Thus $\lambda\mapsto\min_{\pi\in\mathcal{U}}\hat{J}(\pi, \lambda)$ is a strictly concave function, so its stationary point
\[
\lambda^*=\frac{M_2+(P_{0}h_{2,0}^2+M_1)z-P_{0}h_{2,0}(x-h_{1,0})}{1-P_{0}h_{2,0}^2-M_1}
\]
is the unique maximizer, which leads to
\begin{align*}
\max_{\lambda\in\mathbb{R}}\min_{\pi\in\mathcal{U}}\hat{J}(\pi, \lambda)&=\frac{P_{0}h_{2,0}^2+M_1}{1-P_{0}h_{2,0}^2-M_1}z^2+2\frac{M_2-P_{0}h_{2,0}(x-h_{1,0})}{1-P_{0}h_{2,0}^2-M_1}z\\
&\quad +M_3+P_{0}(x-h_{1,0})^2+\frac{[M_2-P_{0}h_{2,0})(x-h_{1,0})]^2}{1-P_{0}h_{2,0}^2-M_1}\\
&=\frac{P_{0}h_{2,0}^2+M_1}{1-P_{0}h_{2,0}^2-M_1}\Big(z-\frac{P_{0}h_{2,0}(x-h_{1,0})-M_2}{P_{0}h_{2,0}^2+M_1}\Big)^2\\
&\quad-\frac{[M_2-P_{0}h_{2,0}(x-h_{1,0})]^2}{P_{0}h_{2,0}^2+M_1}+M_3+P_{0}(x-h_{1,0})^2.
\end{align*}

The above analysis boils down to the following theorem.
\begin{theorem}
\label{efficientth}
The optimal portfolio of problem \eqref{optm} corresponding to $\E (X(T))=z$, as a feedback function of the time $t$, the wealth level $X$, and the market regime $i$, is
\begin{align*}
\pi^*(t, X, i)&=-\frac{1}{P}\left(\sigma\sigma'\right)^{-1}
\Bigg[\Big(P\mu+\sigma\Lambda\Big)(X-h_1-(\lambda^*+z)h_2)+P\sigma\rho-P\sigma\eta\Bigg], 
\end{align*}
where
\begin{align*}
\lambda^*=\frac{M_2+(P_{0}h_{2,0}^2+M_1)z-P_{0}h_{2,0}(x-h_{1,0})}{1-P_{0}h_{2,0}^2-M_1}.
\end{align*}

The mean-variance frontier is
\begin{align*}
\mathrm{Var} (X(T))
&=\frac{P_0h_{2,0}^2+M_1}{1-P_0h_{2,0}^2-M_1}\Big(\E (X(T))-\frac{P_0h_{2,0}(x-h_{1,0})-M_2}{P_0h_{2,0}^2+M_1}\Big)^2\\
&\quad-\frac{[M_2-P_0h_{2,0}(x-h_{1,0})]^2}{P_0h_{2,0}^2+M_1}+M_3+P_0(x-h_{1,0})^2
\end{align*}
with $0<P_0h_{2,0}^2+M_1<1$.
\end{theorem}
\begin{remark}
If there is no liability, i.e. $b(t,i)\equiv0, \ \rho(t,i)\equiv0$, then $h_1(t,i)\equiv0, \ \eta_1(t,i)\equiv0$, $M_2=M_3=0$ and Theorem \ref{efficientth} degenerates to Theorem 5.11 of \cite{HSX}.
\end{remark}

\end{document}